\documentstyle[12pt,amssymb,amscd]{article} 
\author{Dmitry E. Tamarkin}

\newtheorem{Theorem}{THEOREM}[section]

 \newtheorem{Definition}[Theorem]{DEFINITION}
 
 \newtheorem{Proposition}[Theorem]{PROPOSITION}
 \newtheorem{Lemma}[Theorem]{LEMMA}
 
\newtheorem{Corollary}[Theorem]{COROLLARY}
\title{The Deformation Complex of
a d-algebra is a (d+1)-algebra}
\begin{document}
\newcommand{\be}{\begin{equation}}
\newcommand{\ee}{\end{equation}}
\hsize 15cm
\vsize 20cm
\leftmargin=-4cm
\topmargin=1cm
\maketitle
\def\isom {\cong}
\def \dash {\mbox{-}}
\def \endproof {$\bigtriangleup$}
\def \Der{{\rm Der}\;}
\def \Def{{ \rm def}\;}
\def \hom{\underline{\rm Hom}\;} 
\def \ei{d\dash alg_\infty}
\def \Hom{\underline{\rm Hom}}
\section{Introduction}
In this paper we give a purely algebraic proof 
of a variant of the theorem of Kontsevich \cite{K}
stating that the deformation complex of a $d$-algebra
shifted by $d$ is naturally a $(d+1)$-algebra. 

By a $d$-algebra structure on a vector space $V$ over a field
of characteristic zero 
we mean the analogue of Poisson algebra with its bracket
of degree $1-d$. Our definition coincides with the usual
one for $d\geq 2$. Nevertheless, in this paper '1-algebra'  means
'usual Poisson algebra', not 'associative algebra'.

 Kontsevich's definition of $d$-algebra    
\cite{K} is as follows: a $d$-algebra is an algebra 
over the chain operad of the 
operad of little $d$-disks. But in the same paper he proves
that over any field of characteristic zero 
this chain operad is quasi-isomorphic to its homology operad.
Thus, our theorem
is a particular, characteristic zero case of Kontsevich's theorem. 

   Let us outline the idea of the proof. First, let us discuss
the notion of a deformation Lie algebra of a $d$-algebra.
As in the case of associative, Lie, commutative etc. algebras, first, we need 
the definition of a homotopy $d$-algebra, which
can be found with the help of the theory of  Koszul operads \cite{GK}.
The key point is that the operad governing $d$-algebras is Koszul
as it was shown in \cite{GJ}. The definition of the structure
of a homotopy $d$-algebra on a complex $V$ says that it is the same as a 
differential on the cofree $d$-coalgebra 
$Cofree_d(V[-d])$ cogenerated by $V[-d]$.
Structures of  usual $d$-algebra correspond to quadratic differentials
on $Cofree_d(V[-d])$.  Let $X$ be a $d$-algebra. Denote by 
$X^\lor$ the $d$-coalgebra $Cofree_d(V[-d])$ equipped with the
quadratic differential corresponding to the $d$-algebra structure
on $X$. The reader familiar with the definition of homotopy
associative ( Lie, commutative, etc.) algebra will see that the 
definition of homotopy $d$-algebra is similar.  The 
analogue of $X^\lor$ in those cases is called the bar complex of
$X$.  

The deformation complex $\Def(X)$ of a $d$-algebra $X$ is just the differential
graded Lie algebra of derivations of $X^\lor$. This object admits a
more 'geometric' definition in terms of the infinitesimal neighborhood
of the identity in the 'algebraic group' ${\rm Aut} X^\lor$. Let us explain the  meaning of this.
We will use the technical notion of 
 coproartinian cocommutative coalgebra
(Section \ref{coprrt}) which is dual to the notion
of proartinian local commutative algebra.
 Now, note that the tensor product
of a $d$-coalgebra and a cocommutative coalgebra is naturally
a $d$-coalgebra.  
Define a functor $F'$ from the opposite to the
category of coproartinian cocommutative coalgebras to the category of sets by 
setting $F'(a)={\rm Hom}_{d-coalg}(X^\lor\otimes a,X^\lor)$.
One sees that the composition of any two elements from $F'(a)$ is
well defined. Thus, $F'$ is actually a functor taking values in  
the category of monoids.         
Since we need the neighborhood of identity, define a subfunctor $F$ of $F'$ by taking as $F(a)$
only those morphisms $X^\lor\otimes a\to X^\lor$ for which the through map
$X^\lor\to X^\lor\otimes a\to X^\lor$ is identity. Here the first arrow
is induced by the canonical inclusion of the ground field $k$ to $a$.
One sees that $F$ takes values in groups.
It turns out that $F$ viewed as a functor to the category of sets
is representable: there exists a cocommutative coalgebra $A$ such that
$F(a)\cong {\rm Hom}(a,A)$ naturally in $a$. The associative
composition law
$F\times F\to F$ defines an associative map of cocommutative
coalgebras $A\otimes A\to A$ making $A$ a bialgebra. 
One sees that $A\cong U(\Def(X))$, $U$ meaning the universal
enveloping algebra.

 One can modify the above construction.
Let $X,Y$ be $d$-algebras and $\phi:X^\lor\to Y^\lor$ a morphism.
Define a functor
$$
F^\phi_{X,Y}: Coproartalg^{\rm op}\to Sets
$$ by setting
$F^\phi_{X,Y}(a)\subset {\rm Hom}_{d\dash coalg}(X^\lor\otimes a,Y^\lor)$,
where we take only those morphisms for which the through map
$X^\lor\to X^\lor\otimes a\to Y^\lor$ is $\phi$.    
 Again, this functor is representable. Denote the corresponding
coproartinian coalgebra by $\Hom^{\phi}_c(X,Y)$. 
Let
$$
X^\lor\stackrel{\phi}{\to} Y^\lor\stackrel{\psi}{\to}
Z^\lor
$$
be the sequence of morphisms. Then we have a natural composition
\begin{equation}\label{cmp}
\Hom^{\phi}_c(X,Y)\otimes \Hom^{\psi}_c(Y,Z)\to 
\Hom^{\psi\circ\phi}_c(X,Z).
\end{equation}

Now we can take the advantage of the fact that there is a natural tensor product on the category of $d$-coalgebras. This means that the functor
$F^{\phi}_{X,Y}$ can be extended to the category of coproartinian
$d$-coalgebras (the definition remains the same but we allow $a$ to be
any coproartinian $d$-coalgebra). This functor is also
representable. Denote the corresponding coproartinian $d$-coalgebra by
 $\Hom^{\phi}(X,Y)$. One sees that this construction is similar
to the one of internal homomorphisms. The only difference is that
we work in the formal neighborhood of a given morphism. The most interesting case for us is $\Hom^{Id}(X,X)$. 
The analogue of the composition
morphism (\ref{cmp}) provides us with an associative map of $d$-coalgebras $\Hom^{Id}(X,X)\otimes \Hom^{Id}(X,X)\to \Hom^{Id}(X,X)$. We call this structure 
$d$-bialgebra. One sees that as a $d$-coalgebra (the differential is ignored)
$\Hom^{Id}(X,X)$ is isomorphic
to the cofree $d$-coalgebra cogenerated by $\Def(X)$ (viewed as
a graded vector space). The differential and the associative product
look more sofisticated. 

{\bf Remark} Note that the same construction is applicable to
the category of associative algebras. In this case $\Hom^{Id}(X,X)$
is the Hopf algebra isomorphic to the Hopf algebra
on the tensor coalgebra $T(C^{\bullet}(X,X)[1])$ defined in \cite{GJ}
(the so-called $B_\infty$-structure). 

Our next step is to show that the structure of a $d$-bialgebra
on a cofree $d$-coalgebra cogenerated by a complex $K$ 
implies  the structure of $(d+1)$-algebra on $K[-d]$. Contrary to the case of associative algebras this can be done
in an easy purely algebraic way. This result
applied to $\Hom^{Id}(X,X)$  gives the desired $(d+1)$-algebra
structure on $\Def(X)[-d]$.

Here is the content of the sections. We start with the quick
review of $d$-algebras, their homotopy theory, and their deformations. Also we introduce the notion of the coproartinian coalgebra
which is the dual to the notion of proartinian algebra. 
    In Section 2 we construct the  functor $F^{\phi}_{X,Y}$ 
 and show that it is representable. 
The corollary of this is the fact that we have a structure of 
$d$-bialgebra on a cofree $d$-coalgebra cogenerated by the deformation
complex of a $d$-algebra.
In section 3 we prove that this structure implies the desired structure of $(d+1)$-algebra on the shifted deformation complex.

The author would like to thank B. L.  Tsygan and P. Bressler
for their help.

\section{$d$-Algebras, homotopy $d$-algebras, and their deformations}
In this section we will recall the notions outlined in its title.
All the definitions are parallel to the ones for associative, commutative, or
Lie algebras. 
\subsection{$d$-Algebras } \label{dalg}
 A structure of  ${d}$-{\it algebra} on a complex $V$ of vector
spaces over a field $k$ of characteristic zero
is given by 
\begin{enumerate}
\item[1] a dg commutative associative product $\cdot:S^2V\to V$ ;
\item[2] a map 
\begin{equation}\label{commutator}
\{\}:\Lambda^2(V[d-1])\to V[d-1],
\end{equation}
 which turns $V[d-1]$ into DGLA (Differential Graded Lie Algebra).
\end{enumerate}
By  abuse of notation we will also denote by
$\{,\}$  the  degree $1-d$ map $S^2V\to V$ if $d$ is odd
($\Lambda^2V\to V$ if $d$ is even) corresponding to the map
 (\ref{commutator}).

These operations must satisfy
the Leibnitz identity
 $$
\{ab,c\}=a\{b,c\}+(-1)^{|b|(|c|+d-1)}\{a,c\}b.
$$

A $d$-{\it algebra with unit} is a $d$-algebra with a marked  
element $1$ which is the unit with respect to the product and
its bracket with any element vanishes.
The structure of $d$-algebra (resp. $d$-algebra
with unit) is governed by an operad which is denoted
by $e_{d}'$ (resp. $e_{d}$). This means that a $d$-algebra 
(resp. $d$-algebra
with unit) is the same as an algebra 
over the operad $e_d'$ (resp. $e_d$).
As usual, we define a morphism of $d$-algebras (with unit)  $V$ and $W$ as
a morphism of complexes $V\to W$ which respects all the operations (and the units). Therefore, we have the categories
 $d\dash alg$ of $d$-algebras and $d\dash alg_1$ of $d$-algebras with unit.

     \subsection{$d$-Coalgebras} This structure is dual
to the structure of $d$-algebra.
A structure of $d$-coalgebra (with counit) on a complex $V$  is specified by
a dg cocommutative coassociative coproduct $\Delta: V\to S^2V$
(with counit) $\epsilon:V\to k$ and a Lie cobracket $\delta: V[1-d]\to
\Lambda^2V[1-d]$, satisfying the conditions dual to the ones for 
$d$-algebras. $d$-Coalgebras (resp. $d$-coalgebras with counit) form a category $d\dash coalg$
(resp. $d\dash coalg_1$).
      \subsection{Free $d$-algebras} 
Let $Complexes$ be the category of complexes
 and their 
morphisms of degree zero.
We have the forgetful functor
$Obl:d\dash alg_1\to complexes$ which takes a $d$-algebra to its underlying
complex. This functor has the left adjoint. We will denote it
$Free_d$ or simply $Free$. 
The theory of operads says that 
$$
Free_d(V)=\bigoplus\limits_{n=0}^\infty e_d(n)\otimes_{S_n}
V^{\otimes n}.
$$

For a complex $V$, let $FreeLie(V)$ denote the free graded Lie algebra   
generated by $V$ with the differential induced by the one on $V$. 
Denote by $OblFreeLie(V)$ the underlying complex of $FreeLie(V)$. 
Then
on the level of complexes we have the canonical isomorphism
$$
 S^{\bullet}(FreeLie(V[d])[-d])\cong OblFree_d(V).
$$
\subsection{Cofree coproartinian $d$-coalgebras} 
Similarly, we have the forgetful functor
$Obl:
d\dash coalg\to complexes$. Unfortunately,  it does not have an adjoint functor.
The reason is that the linear dual to the free $d$-algebra  generated by even
finite dimensional space $V$ is not a $d$-coalgebra.
This can be cured either by passage to the category of topological 
$d$-coalgebras
or by passage to coproartinian coalgebras. We choose the second option.
\subsubsection{Coproartinian $d$-coalgebras}
\label{coprrt}
 A $d$-coalgebra with counit
$V$ is called ${\em coproartinian}$ if there exists an exhaustive increasing filtration $F$ of $V$ by sub-coalgebras with counit and the following
conditions are satisfied.
\begin{enumerate}
\item[1] $F^0V$ is a one-dimensional coalgebra.
\item[2] Let $1\in F^0V$ be a unique element such that $\epsilon(1)=1$,
  where $\epsilon$ is the counit. Then for any $x\in F^1V$ we have
\be \label{coartcocomm1}
\Delta x=1\otimes x+x\otimes 1;
\ee
\be\label{coartlie1}
 \delta x=0,
\ee
 and for any
$
x\in F^iV, \quad i\geq2
$
 we have 
\be\label{coartcocomm2}
\Delta x-1\otimes x-x\otimes 1\in F^{i-1}V\otimes F^{i-1}V
\ee
and 
\be\label{coartlie2}
\delta x\in F^{i-1}V\otimes F^{i-1}V.
\ee
\end{enumerate}
\begin{Lemma}\label{grlike} 	If $V$ is a coproartinian coalgebra, 
it has the only grouplike element $e$. Therefore,
all  filtrations on $V$ that make it coproartinian
have the same term $F^0V=ke$.
\end{Lemma}

\noindent{\em Proof.} Let $F$ be a filtration on $V$ that
makes it coproartinian. Then we take the element $1\in F^0V$
such that $\epsilon(1)=1$,  as we did it in the definition of
coproartinian coalgebras. It is grouplike. Let $i$ be the minimal
number such that $F^i(V)$ contains a grouplike element $x\neq 1$.
Clearly, $i>0$.
We have from (\ref{coartcocomm2}) that 
$$
(x-1)\otimes (x-1)=\Delta x-1\otimes x -x\otimes 1+1\otimes 1\in 
T^2F^{i-1}(V).
$$ 
Therefore, $x-1\in F^{i-1}V$ and $x\in F^{i-1}V$. Contradiction. 

\begin{Lemma}\label{plus} Let W be a coproartinian coalgebra and let 
$s_n:W\to {\rm Hom}_{S_n}(e_d(n),W^{\otimes n})$ be the structure
maps of $W$ as a $d$-coalgebra. Let
 $$
s_n': W\to {\rm Hom}_{S_n}\big (e_d(n),(W/F^0W)^{\otimes n}\big )
$$
be the composition of $s_n$ with the $n$-th tensor power of the projection $W\to W/F^0W$.
Then $s_n'(F^kW)=0$ whenever $n>k$.
\end{Lemma}

Denote by $d\dash coart$ the full subcategory of $d\dash coalg_1$ formed by 
 coproartinian coalgebras.
 Define the forgetful functor $Obl: d\dash coart\to complexes$
by setting $Obl(V)=V/F^0V$. This is a well-defined functor
by Lemma \ref{grlike}.
\begin{Proposition} The forgetful functor $Obl$ has the right adjoint 

\end{Proposition}
\noindent{\em Proof.} 
Define the right adjojnt functor $Cofree_d$ as
 $$
Cofree_d(V)=\bigoplus_{n=0}^{\infty} {\rm Hom}_{S_n}(e_d(n),V^{\otimes n}).
$$
The structure morphisms of $e_d$ define the canonical $d$-coalgebra structure on
it. Also, we have a filtration on it by the tensor powers of $V$. So, it is a coproartinian coalgebra.
Denote $p: \bigoplus_{n=0}^{\infty} {\rm Hom}_{S_n}(e_d(n),V^{\otimes n})\to V$ the projection onto
$Hom_{S_1}(e_d(1),V)\cong V$.

  Let $W$ be another coproartinian coalgebra. Consider $Hom_{d\dash coart}(W,Cofree_d(V))$. The composition with $p$
defines the map 
$$
p_*:Hom_{d\dash coart}(W,Cofree_d(V))\to
{\rm Hom}_{Complexes}(W,V).
$$
 By Lemma \ref{grlike}
$F^0W$ goes to $F^0Cofree_d(V)$ under any coalgebra morphism. Therefore, the map $p_*$ induces a map
$$
p_*':Hom_{d\dash coart}(W,Cofree_d(V))\to
{\rm Hom}_{Complexes}(W/F^0(W),V),
$$
 and we need to prove that
this map is an isomorphism. For this, we 
define a map 
$$
q:{\rm Hom}_{Complexes}(W/F^0(W),V)\to {\rm Hom}_{d\dash coart}(W,Cofree_d(V)) 
$$ 
and show that it is inverse to $p_*'$. By Lemma \ref{plus}, we have a well defined structure
map $s:W\to \oplus_{n=1}^{\infty} {\rm Hom}_{S_n}(e_d(n),(W/F^0W)^{\otimes n})$. 
Let $r\in {\rm Hom}_{Complexes}(W/F^0(W),V)$.
Define $q'(r)$ as the  composition of $s$ and tensor powers of $r$
and set $q(r)(w)=q'(r)(w)+\epsilon(w)\cdot 1$, where $w\in W$ and $1$ is the grouplike element in $Cofree_d(V)$.
 One sees that $q$ is a coalgebra morphism and that $q$ is inverse to $p_*'$.
\endproof 

\subsubsection{Differentials on $d$-coalgebras} Let $W$ be a $d$-coalgebra. A graded map 
 $D:W\to W$
is called {\em a derivation of} $W$, if it is a derivation of $W$ with respect to its 
cocomutative and Lie coalgebra structures 
(for a moment, we forget the differential on $W$). If the grading  $|D|=1$, and
$D^2$=0, then $D$ is nothing else but a differential on $W$. All derivations of $W$ form a graded Lie algebra
$\Der(W)$ with the commutator $[D,E]=D\circ E-(-1)^{|D||E|}E\circ D$, $\circ$ meaning the composition of maps $W\to W$.
The differential $d$ on $W$ is an element of this Lie algebra, and $[d,d]=0$. Define the differental 
$\delta$ on $\Der(W)$  by $ \delta x=[d,x]$. This turns $\Der(W)$ into a DGLA.
 
\subsubsection{Derivations on cofree coalgebras } Let $Cofree_d(V)$ be the cofree $e_d$-coalgebra
cogenerated by a graded space $V$ with zero differential. Let $p:Cofree_d(V)\to V$ be the canonical projection onto cogenerators.
Let $D\in \Der(Cofree_d(V))$. Then we have the corestriction $cor(D)=p\circ D: Cofree_d(V)\to V$.
\begin{Proposition} The map 
$$
cor:  \Der(Cofree_d(V))\to  {\rm Hom}_k(Cofree_d(V),V)
$$
is an isomorphism of graded vector spaces.
\end{Proposition}
The $Cofree_d(V)$ is graded by the tensor powers of $V$ so that 
$$
gr_i Cofree_d(V)\cong {\rm Hom}_{S_i}(e_d(i),V^{\otimes i}).
$$

Define a subspace $\Der_{1,2}\subset \Der(Cofree_d(V))$  consisting of the elements $x$
such that $cor(x)(gr_i(Cofree_d(V)))=0$ for all $i$ except 1 and 2.
\begin{Proposition}\label{differ} The set of the elements $x\in \Der_{1,2}$ such that $|x|=1$
and $[x,x]=0$ is in $1-1$ correspondence with the structures of differential
$d$-algebras on the shifted graded vector space $V[-d]$. 
\end{Proposition}

\noindent{\em Proof.} Take the components of $cor(x)$   $x_1:V\to V$ and $x_2: {\rm Hom}_{S_2}(e(2), V^{\otimes 2})\to V$.
One sees that $x_1$ is a differential on $V$ and, hence, on $V[-d]$. Recall that $e_d(2)$ is generated by
  the elements  $m,b$  corresponding to the commutative product and the bracket. Therefore,
$j:{\rm Hom}_{S_2}(e(2), V^{\otimes 2})\cong S^2V\oplus \Lambda^2(V[1-d])[d-1]$. Since $x_2$ has degree 1,
it defines under $j$ a map of degree zero $k:\Lambda^2(V[-1])[1]\oplus S^2(V[-d])[d]\to V$. The condition 
$[x,x]=0$
is equivalent to the following:
\begin{enumerate}
\item[1] each of the restrictions of 
$k$ 
$$
b:\Lambda^2(V[-1])[1]\to V \mbox{ and } m:S^2(V[-d])[d]\to V
$$
 are compatible with the differential on
$V$ defined by $x_1$;
\item[2] The map $m$ is the commutative product on $V[-d]$ and $b$ is the Lie bracket on $V[-1]$;
the maps $m,b$ define a structure of $d$-algebra on $V[-d]$.
\end{enumerate}
 Whence the statement of the proposition
\endproof

\begin{Definition} For a $d$-algebra $V$ define the coproartinian $d$-coalgebra  $V^\lor$ as the coalgebra
$Cofree_d(V[d])$ with the differential corresponding to the $d$-algebra structure on $V\cong V[d][-d]$
by Proposition \ref{differ}
\end{Definition}
\subsection{Homotopy $d$-algebras}
\begin{Definition} A structure of homotopy $d$-algebra on a graded vector space $V$ is a differential 
of  the coalgebra $Cofree_d(V[d])$
vanishing on $1\in Cofree_d(V[d])$. 
\end{Definition}
For a homotopy $d$-algebra $V$ we denote by $V^\lor$ the corresponding differential cofree
 coalgebra.
\begin{Definition}
 A morphism of homotopy $d$-algebras is a morphism of the corresponding
differential cofree coalgebras.
\end{Definition}
Thus,  homotopy $d$-algebras form a category. 

We see that any $d$-algebra $V$ defines a homotopy 
$d$-algebra $V^\lor$.
But the set of morphisms between  two $d$-algebras viewed as homotopy $d$-algebras
is wider than the set of usual morphisms between them. In other words, 
 we have an injection
$ {\rm Hom}_{d\dash alg}(V,W)\to {\rm Hom}_{d\dash coart}(V^\lor,W^\lor)$.

 For a homotopy $d$-algebra $V$ the {\em linear part} of the differential $d$ on $V^\lor$
is the restriction of $d$ onto $gr_1(V^\lor)\stackrel{def}{=}gr_1(Cofree_d(V[d]))\cong V[d]$. It 
takes values in $gr_1V^\lor$ and defines a differential on $V[d]$. A structure
of a homotopy $d$-algebra on a complex $V$ is by definition a differential
on $ Cofree_d(V[d])$ such that its linear part coincides with the differential on $V$.

From the operadic point of view, the structure of a homotopy $d$-algebra on a complex
is governed by a dg -operad. Denote it by $he_d$. Let $e_d'$ be the operad governing
$d$-algebras without unit.
The fact that any usual $d$-algebra is also a homotopy $d$-algebra reflects in a map
$p:he_d\to e_d'$
It is known that $he_d$ is a free operad and that $p$ is a quasiisomorphism of operads.
Thus $he_d$ is a free resolution of $e_d'$.

\subsubsection{Deformation Lie algebra}
\begin{Definition} let $V$ be a homotopy $e_d$-algebra. Define its deformation
Lie algebra $\Def(V)=\Der(V^\lor)$ with the differential being the bracket with the differential
on $V^\lor$. 
\end{Definition} 
\section{Infinitesimal Internal Homomorphisms}
In this section first we define the tensor product of 
(coartinian) $e_d$-coalgebras, and then construct  a substitute
for the internal homomorphisms.
\subsection{Tensor product of $e_d$-coalgebras with counit}
Let $V,W$ be $d$-coalgebras with counit. Define  the $d$-coalgebra structure
on $V\otimes W$ as follows. The differential on $V\otimes W$
is the differential on the tensor product of complexes.
The coproduct is defined by $\Delta(v\otimes w)=\epsilon\Delta(v)\otimes \Delta(w)$
and the cocommutator $\delta(v\otimes w)=\epsilon(\delta(v)\otimes \Delta(w))+
\epsilon((-1)^{(d-1)|v|}\Delta(v)\otimes\delta(w))$,
where $\epsilon$ means the sign corresponding to the permutation $(1324)$ of the graded tensor factors:
 $$
\epsilon:V\otimes V\otimes W\otimes W\to V\otimes W
\otimes V\otimes W.
$$ 
 The counit is the tensor
product of counits.
One sees that the tensor product of coproartinian coalgebras
is a coproartinian coalgebra.
\subsection {Internal homomorphisms}
\subsubsection{Useful Lemma}

Let $A$ be a coproartinian $d$-coalgebra; $B$ a 
$d$-algebra without unit. 
We have an injection of $d$-coalgebras
$k\cong F^0A\to A$. Therefore, the factor $A/k$
is naturally a $d$-coalgebra.  
Then $Hom_k(A/k,B)$ is a $d$-algebra,
hence, $Hom_k(A/k,B)[d-1]$ is a Lie algebra. For a Lie algebra $\frak g$ denote $MC(\frak g)=\{x\in \frak g^1:dx+[x,x]/2=0\}$.

\begin{Lemma}\label{useful} There is a natural bijection between
the sets ${\rm Hom}_{d\dash coart}(A, B^\lor)$ and
$MC(Hom_k(A/k,B)[d-1])$.
\end{Lemma}

\noindent{\em Proof.} If we forget about the differentials,
then 
$$
{\rm Hom}_{d\dash coart}(A, B^\lor))\cong Hom_k(A/k,B[d])
\cong  Hom_k(A/k,B)[d-1]^1.
$$
 A direct computation shows that 
the morphisms compatible with the differential correspond
under this identification to the $MC(Hom_k(A/k,B)[d-1])$.
\endproof
\subsubsection{Coalgebra $\hom^\phi(V,W)$}
Let $S$ be a coproartinian dg coalgebra. We have a retraction 
$$
k\cong F^0S\to S\stackrel{\epsilon}{\to} k.
$$
Let $\phi\in {\rm Hom}_{d\dash coart}(V^\lor,W^\lor)$. 
The canonical inclusion $k\to S$ defines a map
$$
h:{\rm Hom}_{d\dash coart}(V^\lor\otimes S,W^\lor)\to
{\rm Hom}_{d\dash coart}(V^\lor,W^\lor).
$$
Set $F^{\phi}_{VW}(S)=h^{-1}(\phi)$. $F^{\phi}_{VW}$ is a 
functor $d\dash coart^0\to Sets$. 
\begin{Proposition} $F^{\phi}_{VW}$ is representable.
\end{Proposition} 

\noindent{\em Proof.} Denote the coalgebra which represents
$F^{\phi}_{VW}$ by  $\hom^\phi(V,W)$. Let us construct it.
Take a $d$-algebra $a={\rm Hom}_k(V^\lor,W)$. By Lemma
\ref{useful} the morphism $\phi$ defines an element
$\phi'\in MC(a[d-1])$
via the inclusion
 ${\rm Hom}(V^\lor/k,W)\to{\rm Hom}(V^\lor,W)$.
 We will denote by the same
letter the corresponding element of degree $d$ in $a$.
Let $a'$ be a $d$-algebra whose operations are the same
as in $a$ but the differential is
$d'x=dx+\{\phi',x\}$, where
$d$ is the differential on $a$. We claim that  
$\hom^\phi(V,W)=a'$. Indeed, Let $b={\rm Hom}_k(V^\lor\otimes S/k,W)$ and $c={\rm Hom}(V^\lor/k,W)$. We have 
$$
{\rm Hom}_{d\dash coart}(V^\lor\otimes S,W^\lor)\cong
MC(b[d-1]).  
$$

We have a retraction 
$$
c\stackrel{G}\to
 b\stackrel{H}\to c
$$ induced by the canonical
retraction $k\to S\to k$. Therefore, we have a semidirect 
sum of $d$-algebras
\begin{equation}\label{split} 
b\cong c + {\rm Hom}_k(S/k,a).
\end{equation}

The  map of Maurer-Cartan elements induced by $H$
is the map $h$. Therefore, $F^{\phi}_{VW}(S)$ can be alternatively described as the set of  $x\in MC(b[d-1])$, $H(x)=\phi'$. Using the splitting (\ref{split}), we write $x=\phi'+s$, $s\in {\rm Hom}_k(S/k,a)$. The Maurer-Cartan equation reads as
$ds+\{\phi',s\}+\{s,s\}/2=0$. This is the same as to say that $s$ viewed as an 
element of ${\rm Hom}_k(S/k,a')$ is  a  morphism of $d$-coalgebras.
\endproof
\begin{Corollary}
There is a natural map 
\begin{equation}\label{compositia}
\circ:\hom^\phi(U,V)\otimes\hom^\psi(V,W)\to 
\hom^{\psi\circ\phi}(U,W).
\end{equation}
This map is associative, meaning that 
the maps $\circ(\circ\otimes Id)$ and $\circ(Id\otimes \circ)$ from
$ \hom^\phi(U,V)\otimes\hom^\psi(V,W)\otimes\hom^\chi(W,X)$ to
$\hom^{\chi\circ\psi\circ\phi}(U,X)$
 coincide.
\end{Corollary}

\noindent{\em Proof.} We have a natural composition map
\begin{eqnarray*}
{\rm Hom}_{d\dash coart}(U^\lor\otimes S, V^\lor)\times 
{\rm Hom}_{d\dash coart}(V^\lor\otimes T, W^\lor)\\
\to {\rm Hom}_{d\dash coart}(U^\lor\otimes S\otimes T, W^\lor).
\end{eqnarray*}

This map induces a morphism of the  functors $d\dash coart^0\times 
d\dash coart^0\to Ens:$   
\begin{equation}\label{morphism}
F^\phi_{U,V}\times F^\psi_{V,W}\to F^{\psi\circ\phi}_{UW}\circ
\bigotimes,
\end{equation}
 where $\bigotimes: d\dash coart\times 
d\dash coart\to d\dash coart$ is the tensor product.
We have the element
$$
Id\in  F^\phi_{U,V}(\hom^\phi(U,V))\cong 
{\rm Hom}(\hom^\phi(U,V),\hom^\phi(U,V)).
$$
Similarly, we have the element $Id\in  F^\psi_{V,W}(\hom^\psi(V,W))$.
Define the morphism (\ref{compositia}) as the image of $Id\times Id$
under the morphism (\ref{morphism}).
\endproof
\section{$\hom^{Id}(V,V)$ and $\Def(V)$.}
\subsection{$\hom^{Id}(V,V)$ is a $d$-bialgebra }
\begin{Definition}
A structure of coproartinian $d$-bialgebra on a complex $V$ is
\begin{enumerate}
\item[1] a structure of a coproartinian $d$-coalgebra on $V$;
\item[2] a morphism of coalgebras $m:V\otimes V\to V$ such that 
it defines an associative product on $V$ with unit being the
grouplike element of $V$.
\end{enumerate}
\end{Definition}
\begin{Proposition} Let $V$ be a $d$-algebra. Then  $\hom^{Id}(V,V)$
is naturally a coproartinian $d$-bialgebra.
\end{Proposition}

{\em Proof.} The product is given by the composition morphism
(\ref{compositia})
\endproof
\subsection{Restriction of $F^{Id}_{VV}$ onto the subcategory
of coproartinian cocommutative coalgebras}
\begin{Definition} A coproartinian cocommutative coalgebra
is an object  of $d\dash coart$ such that its cocommutator is 0.
Denote by $coart$ the corresponding full subcategory
of $d\dash coart$
\end{Definition} 
\begin{Proposition} The inclusion functor $I:coart\to d\dash coart$
has the right adjoint $C$. The coalgebra $C(a)$ is the biggest
cocommutative subcoalgebra of ${\rm Ker}\delta$, where 
$\delta$ is the cocommutator on $a$
\end{Proposition}
Note that for a coproartinian  $d$-bialgebra $X$,
$C(X)$ is naturally a cocommutative Hopf algebra.
\begin{Corollary} The restriction of the functor 
$F^{Id}_{VV}$ to the category of the coalgebras 
is represented by $C(\hom^{Id}(V,V))$. The associative
product on $C(\hom^{Id}(V,V))$ corresponds to the composition
of functors (\ref{morphism}).
\end{Corollary}
On the other hand, since $\Def(V)$ is the Lie algebra of 
the group of  authomorphisms of $V^\lor$, we have
\begin{Proposition} The restriction of the functor 
$F^{Id}_{VV}$ to the category of the coalgebras 
is presented by the universal enveloping algebra
$U(\Def(V))$ viewed as a cocommutative coalgebra.
The associative product on $U(\Def(V))$ corresponds to 
the morphism of functors (\ref{morphism}).
\end{Proposition}
Thus,
\begin{Corollary} We have a canonical isomorphism
of Hopf algebras $C(\hom^{Id}(V,V))\to U(\Def(V))$
\end{Corollary}

Let us express $C(\hom^{Id}(V,V))$ explicitly.
Recall that as a $d$-coalgebra  
$\hom^{Id}(V,V)\cong a'^\lor$, where $a'$ is the $d$-algebra
${\rm Hom}_k(V^\lor,V)$ with the differential twisted by
 the corestriction of $Id$: $V^\lor\stackrel{Id}\to V^\lor
\stackrel{cor}\to V$ which is, of course,
just the map $cor\in a'$. One sees that
for a complex $X$, $C(Cofree_d(X))
\cong
S(X)$, where $S(X)$ is a cofree coproartinian cocommutative
coalgebra cogenerated by $X$. Therefore, $
C(\hom^{Id}(V,V))\cong S(a'[d])$. One sees that as a complex
$a'[d]$ is isomorphic to $\Def(V)$ and that the composition
$S(a'[d])\isom C(\hom^{Id}(V,V))\to U(\Def(V))$ is the 
Poincar\'e-Birkhoff-Witt isomorphism.

\section{A homotopy $(d+1)$-algebra structure on $\Def(V,V)$}
We can summarize our findings in the following way.
We have a $d$-algebra $a'$, which as a complex is isomorphic to
$\Def(V)[-d]$. Also we have the associative product $a'^\lor\otimes
a'^\lor\to a'^\lor$ which turns it into $d$-bialgebra.  Also we know
 that
the restriction of this product onto 
$S(a'[d])=C(a'^\lor)\subset a'^\lor$
turns $S(a'[d])$ into a Hopf algebra  and
\begin{equation}\label{Hopf}
S(a'[d])\cong U(\Def(V)).
\end{equation}
Let us investigate these structures. First, note that $a'^\lor$
together with the cocommutative coproduct and the associative product is a cocommutative cofree Hopf algebra. It is well known 
that any such an algebra is isomorphic to $U(\frak g)$
for a certain Lie algebra $\frak g$. This Lie algebra
is formed by the primitive elements of $a'^\lor$, its commutator
is the commutator with respect to the associative product and
its differential is the restriction of the differential
on $a'^\lor$.
One sees that for any space $X$ the set of primitive elements of $Cofree_d(X)$
is isomorphic to $CofreeLie(X[1-d])[d-1]$, where $CofreeLie$
means 'the cofree Lie coalgebra cogenerated by'. Thus, we have 
a Lie algebra structure on the  Lie coalgebra
 $CofreeLie(a'[1])[d-1]$. The 
fact that $S(a'[d])$ is a Hopf algebra translates into the fact
that $CofreeLie(a'[1])[d-1]$ is a differential Lie bialgebra, 
meaning that all operations are compatible with the differential
and that the cocycle condition is fulfilled:
\begin{equation}\label{cocycle}
\delta([x,y])=[\delta x,y]+(-1)^{|x|(d-1)}[x,\delta y].
\end{equation}
Note that the only difference between our bialgebra and 
usual bialgebras is that the coproduct has
a nonzero grading. The isomorphism (\ref{Hopf})
means that the restriction of the commutator
on the primitive elements 
\begin{equation}\label{ism}
\Def(V)\cong a'[d]\subset
CofreeLie(a'[1])[d-1]
\end{equation}
 coincides  with the commutator on
$\Def(V)$. Thus, (\ref{ism}) is a morphism of DGLA.

 Now take the chain complex 
$S(CofreeLie(a'[1])[d])$ of
$CofreeLie(a'[1])[d-1]$ as a Lie algebra.

 Note that we have an isomorphism of graded spaces
$S(CofreeLie(a'[1])[d])\cong Cofree_{d+1}(a'[1+d])$.
The cocycle condition (\ref{cocycle})
means
that  the differential on  $ S(CofreeLie(a'[1])[d])$
is compatible with the coproduct and the cocommutator
on $Cofree_{d+1}(a'[1+d])$. This means that we have
a structure of homotopy $(d+1)$-algebra on 
a complex $a'\cong 
\Def(V)[-d]$.  The morphism (\ref{ism}) means that the 
homotopy Lie algebra structure on $\Def(V)$
induced from the homotopy
$(d+1)$-structure  on $\Def(V)[-d]$ coincides with the 
Lie
algebra structure on $\Def(V)$ as the deformation Lie algebra.
As for the commutative binary operation on $\Def(V)$, it is the same as the commutative product on $a'$ as a $d$-algebra. One checks
also that the commutator of $a'$ as a $d$-algebra is homotopy
trivial.

\end{document}